# Testing for changes in polynomial regression

ALEXANDER AUE[1], LAJOS HORVÁTH[2], MARIE HUŠKOVÁ[3] and
PIOTR KOKOSZKA[4]

[1]*Department of Statistics, University of California, One Shields Avenue, Davis, CA 95616, USA.
E-mail: alexaue@wald.ucdavis.edu*

[2]*Department of Mathematics, University of Utah, 155 South 1440 East, Salt Lake City, UT
84112–0090, USA. E-mail: horvath@math.utah.edu*

[3]*Department of Statistics, Charles University, Sokolovská 83, CZ–18600 Praha, Czech Republic.
E-mail: huskova@karlin.mff.cuni.cz*

[4]*Department of Mathematics and Statistics, Utah State University, 3900 Old Main Hill, Logan,
UT 84322–3900, USA. E-mail: Piotr.Kokoszka@usu.edu*

We consider a nonlinear polynomial regression model in which we wish to test the null hypothesis
of structural stability in the regression parameters against the alternative of a break at an
unknown time. We derive the extreme value distribution of a maximum-type test statistic which
is asymptotically equivalent to the maximally selected likelihood ratio. The resulting test is easy
to apply and has good size and power, even in small samples.

polynomials; linear regression; polynomial regression

## 1. Introduction

Testing for structural changes has become a focus of attention in statistics and econometrics, reflected in the broadening of possible settings under consideration and the multitude of test statistics developed to investigate them. For recent references, see Dufour and Ghysels (1996) and Banerjee and Urga (2005). In the present paper, we focus on the (nonlinear) polynomial regression model

$$y_i = \mathbf{x}_i^T \boldsymbol{\beta}_i + \varepsilon_i, \qquad i = 1, \ldots, n, \tag{1.1}$$

where $\{\varepsilon_i\}$ is an error sequence that will be specified below, $\{\boldsymbol{\beta}_i\}$ are $(p+1)$-dimensional deterministic vectors and

$$\mathbf{x}_i = (1, i/n, \ldots, (i/n)^p)^T.$$

In this setting, we are interested in testing the null hypothesis of structural stability against the alternative of a regime switch at an unknown time, that is,

$$H_0: \boldsymbol{\beta}_i = \boldsymbol{\beta}_0, \qquad i = 1, \ldots, n;$$
$$H_A: \text{There is a } k^* \geq 1 \text{ such that}$$

$$\boldsymbol{\beta}_i = \boldsymbol{\beta}_0, \qquad i = 1, \ldots, k^*,$$
$$\boldsymbol{\beta}_i = \boldsymbol{\beta}_A, \qquad i = k^* + 1, \ldots, n,$$

with $\boldsymbol{\beta}_0 \neq \boldsymbol{\beta}_A$. The parameters before and after the change, $\boldsymbol{\beta}_0$ and $\boldsymbol{\beta}_A$, as well as the change-point, $k^*$, are unknown.

Pioneering steps to analyze the structural stability for linear models were carried out by Quandt (1958, 1960). His contributions were subsequently refined and extended by, for example, Hawkins (1989), Andrews (1993), Horváth (1995) and Bai *et al.* (1998). For comprehensive reviews, we refer to Csörgő and Horváth (1997) and Perron (2006). While there is an extensive body of literature available for the asymptotic theory in linear models, one of the key assumptions in this setting requires the regressors $\{\mathbf{x}_i\}$ to be stationary or to satisfy regularity conditions such as

$$\frac{1}{n}\sum_{i=1}^{\lfloor rn \rfloor} \mathbf{x}_i \mathbf{x}_i^T \to \mathbf{C}r, \qquad \text{for all } r \in (0,1), \text{ as } n \to \infty, \tag{1.2}$$

where $\mathbf{C}$ is a $(p+1) \times (p+1)$-dimensional positive definite matrix and $\lfloor \cdot \rfloor$ denotes the integer part. Thus, limit results obtained for linear models do not directly apply to the polynomial regression addressed here.

Groundbreaking work for polynomial regression models is due to MacNeill (1978), who derived, under $H_0$, the asymptotic distribution of partial sums obtained from regression residuals. In the same paper, he also proposed two test statistics to distinguish between stability and regime switches, and studied linear and quadratic functionals of residual partial sums. MacNeill's (1978) results assume independent, identically distributed errors $\{\varepsilon_i\}$. For further work in this setting, see Jandhyala (1993) and Jandhyala and MacNeill (1989, 1997). The assumptions on the error sequence were relaxed by Kuang (1998) who requires $\{\varepsilon_i\}$ to obey a functional central limit theorem. He introduced a generalized fluctuation test based on comparisons of recursive estimates with benchmark estimates by mimicking the test procedure of Kuang and Hornik (1995). The new test statistic is compared to MacNeill's (1978) versions in an extensive simulation study. Hušková and Picek (2005) considered bootstrap methods for more general settings, which include polynomial regression models as a special case. Here, the regresssion function is assumed to be smooth. Our work, however, is more akin to the contributions of Jarušková (1998, 1999) and Albin and Jarušková (2003). In the case $p = 1$, these papers obtain extreme value asymptotics for a test statistic designed to find changes in the regression line by explicitly computing the design matrices appearing in (1.2).

The paper is organized as follows. In Section 2, we precisely state assumptions on the innovations $\{\varepsilon_i\}$, introduce a maximum-type test statistic, formulate the main extreme value asymptotic and motivate the proofs which are given in Sections 4 and 5. Section 3 is devoted to the practical application of our theory.

## 2. Model assumptions and results

In what follows, we study the polynomial regression model as specified in (1.1), in which we wish to distinguish between the structural stability null hypothesis $H_0$ and the breakpoint alternative $H_A$. Assuming for the moment that the error sequence consists of independent, identically distributed normal variables with the same known variance $\sigma^2$, and that the time of change is $k^* = k$, twice the logarithm of the likelihood ratio for this two-sample problem is given by

$$\ell_n(k) = n\sigma^{-2}(\hat{\sigma}_n^2 - [\hat{\sigma}_{k,1}^2 + \hat{\sigma}_{k,2}^2]), \tag{2.1}$$

where

$$n\hat{\sigma}_{k,1}^2 = \sum_{i=1}^{k}(y_i - \mathbf{x}_i^T\hat{\boldsymbol{\beta}}_k)^2 \quad \text{and} \quad n\hat{\sigma}_{k,2}^2 = \sum_{i=k+1}^{n}(y_i - \mathbf{x}_i^T\hat{\boldsymbol{\beta}}_k^*)^2$$

denote the sum of the squared residuals for the first, respectively, second sample, and

$$n\hat{\sigma}_n^2 = \sum_{i=1}^{n}(y_i - \mathbf{x}_i^T\hat{\boldsymbol{\beta}}_n)^2$$

the sum of the squared residuals for the whole sample. Therein, $\hat{\boldsymbol{\beta}}_k$ and $\hat{\boldsymbol{\beta}}_k^*$ are the least squares estimators for $\boldsymbol{\beta}$ based on the first $k$ and the last $n-k$ observations. Elementary algebra implies that, under $H_0$, we have

$$\hat{\sigma}_n^2 - [\hat{\sigma}_{k,1}^2 + \hat{\sigma}_{k,2}^2] = \frac{1}{n}\mathbf{S}_k^T\mathbf{C}_k^{-1}\mathbf{C}_n\tilde{\mathbf{C}}_k^{-1}\mathbf{S}_k,$$

where

$$\mathbf{C}_k = \sum_{i=1}^{k}\mathbf{x}_i\mathbf{x}_i^T, \qquad \tilde{\mathbf{C}}_k = \sum_{i=k+1}^{n}\mathbf{x}_i\mathbf{x}_i^T$$

and

$$\mathbf{S}_k = \sum_{i=1}^{k}\mathbf{x}_iy_i - \mathbf{C}_k\mathbf{C}_n^{-1}\sum_{i=1}^{n}\mathbf{x}_iy_i = \sum_{i=1}^{k}\mathbf{x}_i(y_i - \mathbf{x}_i^T\hat{\boldsymbol{\beta}}_n).$$

Since, in general, the time of change $k^*$ is unknown, we reject the null hypothesis $H_0$ for large values of

$$T_n = \frac{1}{\sigma^2}\max_{p<k<n-p}\mathbf{S}_k^T\mathbf{C}_k^{-1}\mathbf{C}_n\tilde{\mathbf{C}}_k^{-1}\mathbf{S}_k. \tag{2.2}$$

Note that, as a straightforward computation shows, the value of $T_n$ does not change if $\mathbf{x}_i = (1, i/n, \ldots, (i/n)^p)^T$ is replaced by $\mathbf{x}_i^* = (1, i, \ldots, i^p)^T$.

Still assuming that the errors are independent, identically normal, but also that the common variance $\sigma^2$ is unknown, the resulting likelihood ratio is

$$\hat{\ell}_n(k) = (\hat{\sigma}_n^{-2}[\hat{\sigma}_{k,1}^2 + \hat{\sigma}_{k,2}^2])^{n/2},$$

in the case where the time of change is exactly $k = k^*$ (see Andrews (1993) and Csörgő and Horváth (1997), Section 3.1.1). If $k^*$ is unknown, the stability of the regression coefficients is rejected for large values of

$$\hat{T}_n = \max_{p < k < n-p} [-2 \log \hat{\ell}_n(k)]. \tag{2.3}$$

Hansen (2000) studied the asymptotic distribution of the trimmed version

$$T_{n,\delta} = \frac{1}{\hat{\sigma}_n^2} \max_{\lfloor n\delta \rfloor \le k \le n - \lfloor n\delta \rfloor} \mathbf{S}_k^T \mathbf{C}_k^{-1} \mathbf{C}_n \tilde{\mathbf{C}}_k^{-1} \mathbf{S}_k$$

with some $0 < \delta < 1$. The limit distribution of $T_{n,\delta}$ is the supremum of quadratic forms of $(p+1)$-dimensional Gaussian processes with a complicated covariance structure. To state the limit theorem for the truncated statistic $T_{n,\delta}$, we introduce the matrices

$$\mathbf{C}(t) = \left( \int_0^t x^{i+j} \, \mathrm{d}x : 0 \le i, j \le p \right)$$

and

$$\tilde{\mathbf{C}}(t) = \left( \int_t^1 x^{i+j} \, \mathrm{d}x : 0 \le i, j \le p \right).$$

Furthermore, for $t \ge 0$, let $\mathbf{\Gamma}(t) = (\int_0^t x^i \, \mathrm{d}W(x) : 0 \le i \le p)$, where $\{W(t) : t \ge 0\}$ denotes a standard Brownian motion. The proof of the following theorem is due to Hansen (2000).

**Theorem 2.1.** *If* (1.1) *and* (2.2)–(2.6) *hold, then, under* $H_0$,

$$T_{n,\delta} \xrightarrow{\mathcal{D}} \sup_{\delta \le t \le 1-\delta} \mathbf{\Delta}^T(t) \mathbf{C}^{-1}(t) \mathbf{C}(1) \tilde{\mathbf{C}}(t) \mathbf{\Delta}(t)$$

*for all* $\delta \in (0, 1/2)$, *where* $\mathbf{\Delta}(t) = \mathbf{\Gamma}(t) - \mathbf{C}(t)\mathbf{C}^{-1}(1)\mathbf{\Gamma}(1)$ *and* $t \ge 0$.

The use of the truncated statistic $T_{n,\delta}$ requires choosing $\delta$ and so a priori excludes changes close to the end-points of the sample. More importantly, the computation of the asymptotic critical values for $T_{n,\delta}$ is a non-trivial numerical task. A practical use of Theorem 2.1 would require tables for a range of values of $\delta$ and $p$, and we are not aware of any such tables.

In contrast, the definition of $T_n$ allows, at least in principle, to detect changes anywhere in the sample, and asymptotic critical values can be easily computed using Theorem 2.2 below. The theory needed to establish Theorem 2.2 is, however, far from trivial. We are following previous work of Jarušková (1998, 1999) and Albin and Jarušková (2003). These contributions utilize the theory of high level exceedence probabilities for stationary Gaussian processes developed in Albin (1990, 1992, 2001) to obtain extreme value asymptotics. For classical results in the field, see Leadbetter *et al.* (1983); for a more recent survey on extreme value theory with applications in telecommunication and the environment, see Finkenstaedt and Rootzén (2003).

While, for motivational reasons, all test statistics in this section have been introduced and explained for normal error sequences, the resulting test procedures are sensitive with respect to changes in the regression parameters in a much more general setting. We assume that the errors have constant variance and are uncorrelated, that is,

$$\mathrm{E}\varepsilon_i = 0, \qquad \mathrm{E}\varepsilon_i^2 = \sigma^2 \quad \text{and} \quad \mathrm{E}\varepsilon_i\varepsilon_j = 0 \qquad (i \neq j), \tag{2.4}$$

and that there are two independent standard Brownian motions (standard Wiener processes) $\{W_{1,n}(s): s \geq 0\}$ and $\{W_{2,n}(s): s \geq 0\}$ such that

$$\max_{1 \leq k \leq n/2} \frac{1}{k^{1/2-\Delta}} \left| \sum_{i=1}^{k} \varepsilon_i - \sigma W_{1,n}(k) \right| = \mathcal{O}_P(1) \qquad (n \to \infty) \tag{2.5}$$

and

$$\max_{n/2 < k < n} \frac{1}{(n-k)^{1/2-\Delta}} \left| \sum_{i=k+1}^{n} \varepsilon_i - \sigma W_{2,n}(n-k) \right| = \mathcal{O}_P(1) \qquad (n \to \infty) \tag{2.6}$$

with some $\Delta > 0$.

Sequences $\{\varepsilon_i\}$ satisfying the invariance principles (2.5) and (2.6) include many weakly dependent processes including GARCH-type sequences (Aue *et al.* (2006)), mixing sequences (Shao (1993)) and martingale differences (Eberlein (1986)). Note, also, that the assumptions on the error sequence could be further relaxed along the lines of Qu and Perron (2007).

Let $\Gamma(t) = \int_0^\infty \mathrm{e}^{-y} y^{t-1} \,\mathrm{d}y$ denote the Gamma function. The main results are the following limit theorems which establish the asymptotic behavior of $T_n$, $\hat{T}_n$ and $T_{n,\delta}$.

**Theorem 2.2.** *If (1.1) and (2.4)–(2.6) hold, then, under $H_0$, the statistic $T_n$ (2.2) satisfies, for all $x$,*

$$\lim_{n \to \infty} P\left\{ T_n \leq x + 2\log\log n + (p+1)\log\log\log n - 2\log\left(\frac{2^{(p+1)/2}\Gamma((p+1)/2)}{p+1}\right) \right\}$$
$$= \exp(-2\mathrm{e}^{-x/2}).$$

The following corollary to Theorem 2.2 is useful in its practical application.

**Corollary 2.1.** *If (1.1) and (2.4)–(2.6) hold, then, under $H_0$, the statistic $\hat{T}_n$ (2.3) satisfies, for all $x$,*

$$\lim_{n\to\infty} P\left\{\hat{T}_n \leq x + 2\log\log h(n) + (p+1)\log\log\log h(n) - 2\log\left(\frac{2^{(p+1)/2}\Gamma((p+1)/2)}{p+1}\right)\right\}$$
$$= \exp(-2\mathrm{e}^{-x/2}),$$

*where $h(n) = n(\log n)^\gamma$ with an arbitrary real $\gamma$.*

Instead of the statistics $\hat{T}_n$ (2.3), the following statistics can be used:

$$T_{n,1} = \frac{1}{\hat{\sigma}_n^2} \max_{p<k<n-p} \mathbf{S}_k^T \mathbf{C}_k^{-1} \mathbf{C}_n \tilde{\mathbf{C}}_k^{-1} \mathbf{S}_k;$$

$$T_{n,2} = \max_{p<k<n-p} \frac{1}{\hat{\sigma}_{k,1}^2 + \hat{\sigma}_{k,2}^2} \mathbf{S}_k^T \mathbf{C}_k^{-1} \mathbf{C}_n \tilde{\mathbf{C}}_k^{-1} \mathbf{S}_k;$$

$$T_{n,3} = \left(\min_{p<k<n-p}[\hat{\sigma}_{k,1}^2 + \hat{\sigma}_{k,2}^2]\right)^{-1} \max_{p<k<n-p} \mathbf{S}_k^T \mathbf{C}_k^{-1} \mathbf{C}_n \tilde{\mathbf{C}}_k^{-1} \mathbf{S}_k.$$

Indeed, using the Taylor expansion of $\log(1+x)$, $-1 < x < 1$, it is easily verified that, under $H_0$,

$$|\hat{T}_n - T_{n,i}| = \mathcal{O}_P\left(\frac{1}{n}\right) T_{n,i}^2, \qquad i = 1, 2, 3.$$

The statistics $T_{n,1}$, $T_{n,2}$ and $T_{n,3}$ differ only in the built-in estimation of the variance parameter $\sigma^2$. Under the assumptions of Theorem 2.2 and further mild regularity conditions,

$$|\hat{\sigma}_n^2 - \sigma^2| = o_P\left(\frac{1}{\log\log n}\right)$$

and

$$\max_{p<k<n-p} |\hat{\sigma}_{k,1}^2 + \hat{\sigma}_{k,2}^2 - \sigma^2| = o_P\left(\frac{1}{\log\log n}\right)$$

as $n \to \infty$. Hence, the results of Theorem 2.1, Theorem 2.2 and Corollary 2.1 also apply to the test statistics $T_{n,1}$, $T_{n,2}$ and $T_{n,3}$.

The proof of Theorem 2.2 and its corollary is prepared in Section 4 and completed in Section 5. In Section 4, we show that the error sequence can be replaced with independent, identically distributed normal random variables. It will be shown that the limit distribution of $T_n$ is related to the maximum of the sum of squares obtained from continuous-time Gaussian processes which are integrals of Legendre polynomials with respect to the same Brownian motion. The extreme value asymptotics can then be derived utilizing a result of Aue *et al.* (2007). Next, however, we discuss the practical application of our theory.

## 3. Application to finite samples

The goal of this section is twofold: to explain the application of the test procedure and to give an idea of its performance in finite samples.

By Corollary 2.1, the test rejects the null hypothesis at level $\alpha$ if $\hat{T}_n > c(n, \alpha)$, where the critical value $c(n, \alpha)$ is computed via

$$c(n, \alpha) = -2\log(-0.5\log(1-\alpha)) + g(n, p, \gamma),$$

where, setting $h(n) = n(\log n)^\gamma$,

$$g(n, p, \gamma) = 2\log\log h(n) + (p+1)\log\log\log h(n) - 2\log\left(\frac{2^{(p+1)/2}\Gamma((p+1)/2)}{p+1}\right).$$

The symbol log stands for the *natural* logarithm, $p$ is the order of the polynomial regression to be fitted and $\gamma$ is a real constant used to calibrate the size (see below).

The statistic $\hat{T}_n$ is easy to compute in any statistical software package. Denote by $\hat{s}_n$, $\hat{s}_{k,1}$ and $\hat{s}_{k,2}$ the residual standard deviations obtained from fitting a polynomial regression to, respectively, the whole data set, the first $k$ observations and the last $n-k$ observations. Then, after some simple algebra, we obtain

$$\hat{T}_n = -n\left[\min_{p<k<n-p}\{\log((k-p)\hat{s}_{k,1}^2 + (n-k-p)\hat{s}_{k,2}^2)\} - \log(n-p) - 2\log(\hat{s}_n)\right]. \quad (3.1)$$

We recommend reducing the range over which the minimum is taken by one data point on each side, that is, to use the minimum over $k = p+2, \ldots, n-p-2$. An implementation in R (for $p = 2$) is displayed in Figure 1.

The test described above has good size, even in small samples. We considered the linear regression $y_i = \beta_0 + \beta_1 x_i + \varepsilon_i$ and the quadratic regression $y_i = \beta_0 + \beta_1 x_i + \beta_2 x_i^2 + \varepsilon_i$ with i.i.d. standard normal errors $\varepsilon_i$. In both cases, we set $\beta_0 = 1$. For the linear regression, we varied the slope $\beta_1$. For the quadratic regression, we set $\beta_1 = 0$ and varied the coefficient $\beta_2$ of the quadratic term. The empirical sizes based on one thousand replications are displayed in Table 1 – they are within two standard errors of the nominal sizes. Note that for the linear regression we used $\gamma = 0$, as suggested by our main Theorem 2.2, but for the quadratic regression, we used $\gamma = 1$. If $\gamma = 0$ is used for $p = 2$, the rejection rates exceed the nominal size by a few percentage points for the samples sizes we considered.

To assess the power, we used the same general setting as for the size study and considered the following change-point models:

$$\begin{aligned}
p = 1 : \beta_0 &= 1, & \beta_1 &= 1 & &\text{for } i \leq k^*;\\
\beta_0 &= 0, & \beta_1 &= 0 & &\text{for } i > k^*;\\
p = 2 : \beta_0 &= 1, & \beta_1 &= 0, & \beta_2 &= 2 & &\text{for } i \leq k^*;\\
\beta_0 &= 0, & \beta_1 &= 0, & \beta_2 &= 0 & &\text{for } i > k^*.
\end{aligned}$$

**Table 1.** Empirical size (in percent) of the test derived from Corollary 2.1

|   | $p=1$ ($\gamma=0$) | | | | | | $p=2$ ($\gamma=1$) | | | | | |
|---|---|---|---|---|---|---|---|---|---|---|---|---|
|   | $\beta_1=0.5$ | | $\beta_1=1.0$ | | $\beta_1=2.0$ | | $\beta_2=0.5$ | | $\beta_2=1.0$ | | $\beta_2=2.0$ | |
| $n$ | 10% | 5% | 10% | 5% | 10% | 5% | 10% | 5% | 10% | 5% | 10% | 5% |
| 50 | 10.3 | 5.3 | 10.2 | 6.4 | 09.5 | 5.4 | 10.9 | 6.7 | 11.7 | 7.0 | 09.4 | 5.1 |
| 100 | 10.2 | 5.2 | 10.0 | 5.7 | 11.4 | 6.4 | 09.5 | 5.9 | 08.7 | 5.8 | 10.0 | 5.9 |
| 200 | 11.7 | 6.3 | 11.4 | 6.0 | 09.3 | 6.2 | 10.9 | 5.3 | 10.2 | 5.6 | 10.4 | 6.6 |

We considered $k^* = n/2$ and $k^* = n/5$.

Even though these changes appear large, one must keep in mind that they are, in fact, small relative to the standard deviation of the errors. Examples of scatterplots of the simulated data are shown in Figures 2 and 3. It is seen that in many cases, a change-point is not apparent without prior information on its existence. For example, in Figure 2, $n = 200, k^* = 40$, the data appears to fit a linear regression with a negative slope. A visual examination of the remaining panels reveals that the changes we consider are difficult to identify by eye. Nevertheless, Table 2 shows that these changes can be detected with non-trivial power. As for most change-point detection procedures, it is easier to detect a change in the middle of the sample.

In summary, the results of this section show that our procedure is very easy to use in practice and gives good results, even in small samples.

```
hT2=function(Y) #p=2, Y - responses
{
n=length(Y); X1=(1:n)/n; X2= X1*X1
v=rep(-1, n)
for(k in 4:(n-4) )
    {
      lmk1=lm(Y[1:k]~X1[1:k]+X2[1:k])
      lmk2=lm(Y[(k+1):n]~X1[(k+1):n]+X2[(k+1):n])
     sk1=summary(lmk1)$sigma
     sk2=summary(lmk2)$sigma
     v[k]=(k-2)*sk1*sk1+(n-k-2)*sk2*sk2
    }
v=v[4:(n-4)]
lmA=lm(Y X1+X2); sA=summary(lmA)$sigma
-n*(min(log(v)) - log(n-2) - 2*log(sA))
}
```

**Figure 1.** R code for computing the statistic $\hat{T}_n$ (2.3) for $p=2$.

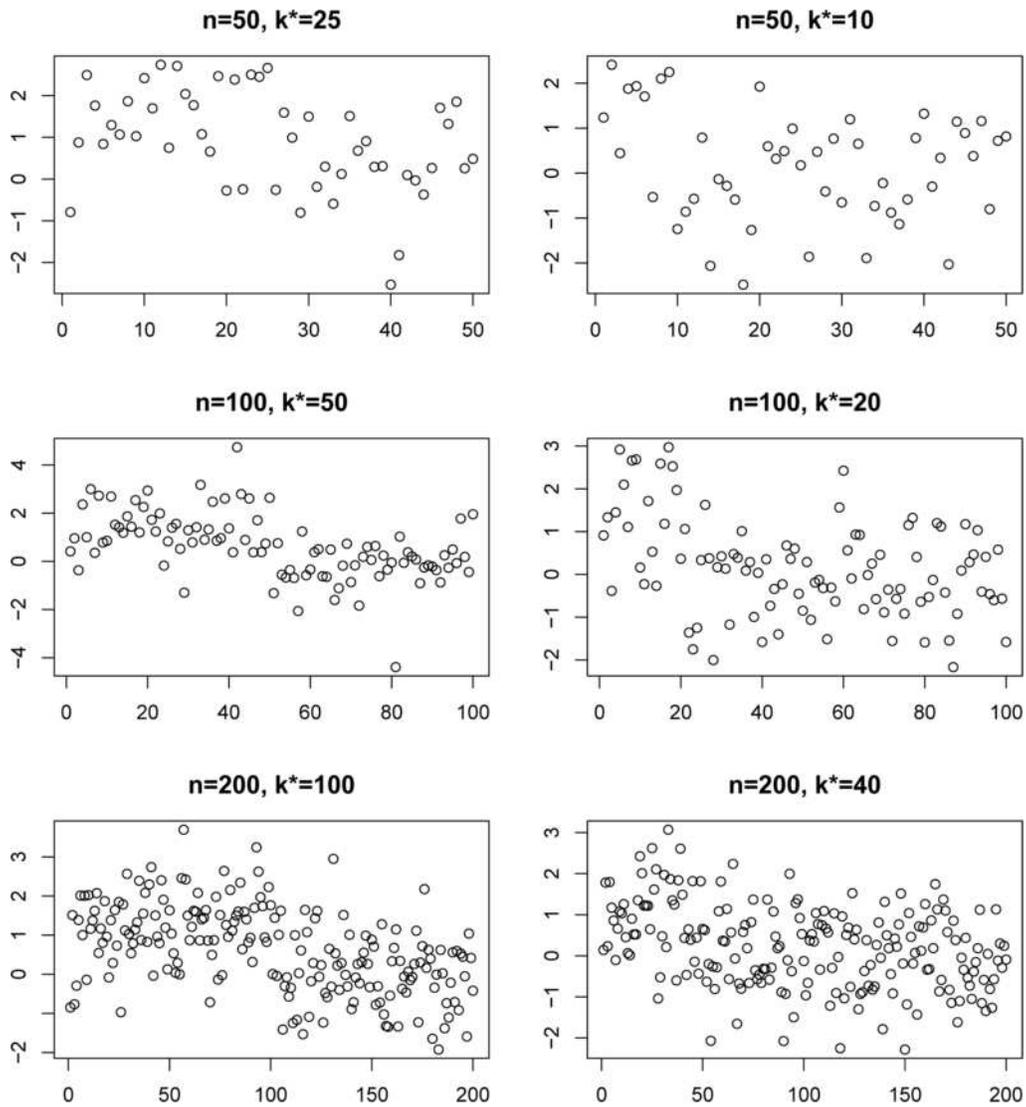

**Figure 2.** Scatterplots for linear regressions with change-points.

## 4. Asymptotic representation of the test statistic

The aim of this section is to derive a more convenient version of the test statistic $T_n$, which includes (a) a transformation that uses only the range of those time-lags $k$ contributing to the extreme value asymptotic, and (b) a continuous-time modification which involves the integration of normalized Legendre polynomials with respect to standard Brownian motion. Throughout, we work under the null hypothesis $H_0$.

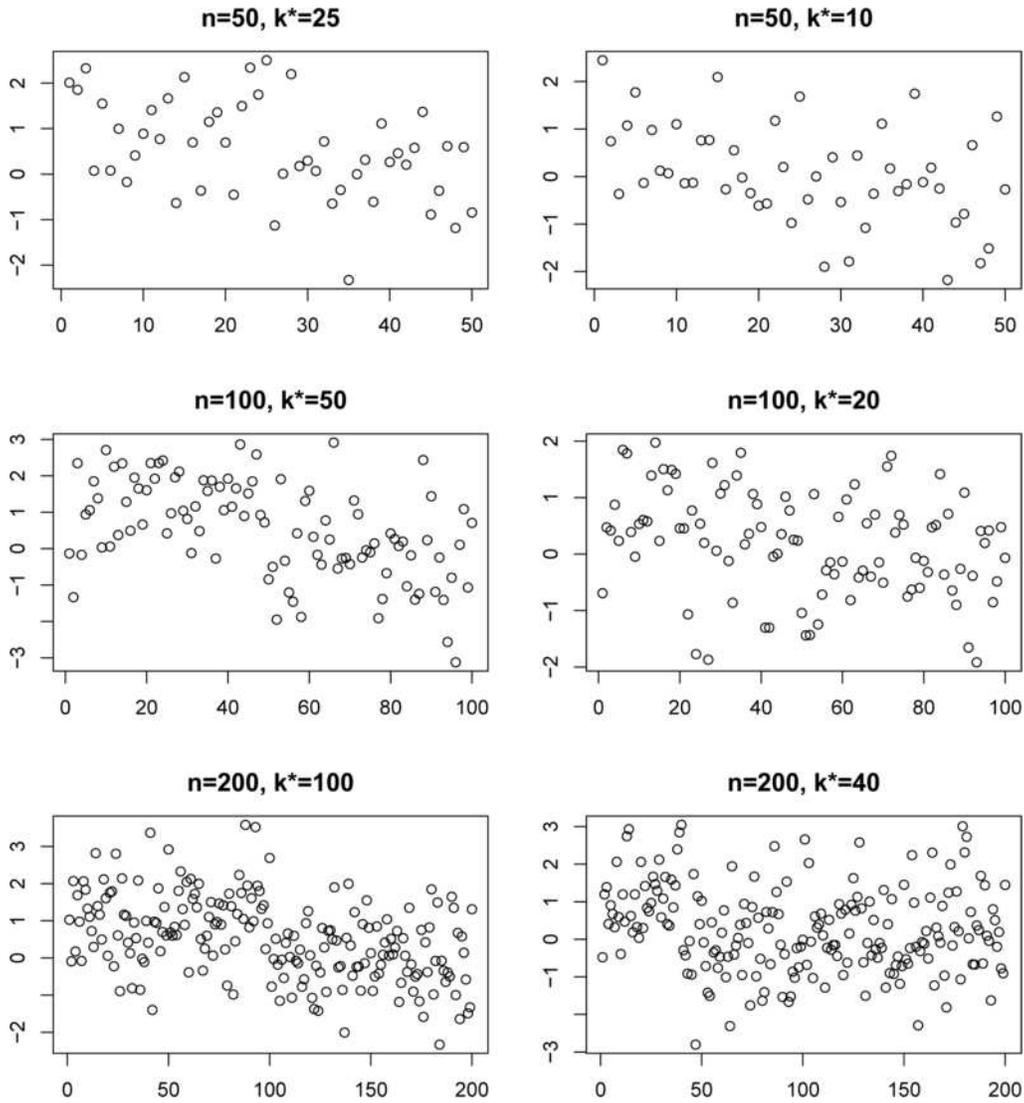

**Figure 3.** Scatterplots for quadratic regressions with change-points.

The following auxiliary lemma will prove useful. It is a consequence of the approximations assumed in (2.5) and (2.6).

**Table 2.** Empirical power (in percent) of the test derived from Corollary 2.1

|   | $p=1$ $(\gamma=0)$ | | | | $p=2$ $(\gamma=1)$ | | | |
|---|---|---|---|---|---|---|---|---|
|   | $k^*=n/2$ | | $k^*=n/5$ | | $k^*=n/2$ | | $k^*=n/5$ | |
| $n$ | 10% | 5% | 10% | 5% | 10% | 5% | 10% | 5% |
| 50  | 45.4 | 34.1 | 40.1 | 29.0 | 36.5 | 28.5 | 16.6 | 10.2 |
| 100 | 80.5 | 71.6 | 67.1 | 56.6 | 68.0 | 57.9 | 24.2 | 17.1 |
| 200 | 99.2 | 98.7 | 94.3 | 91.4 | 96.1 | 94.2 | 46.4 | 35.1 |
| 400 | 100  | 100  | 100  | 98.2 | 100  | 100  | 78.7 | 71.0 |

**Lemma 4.1.** *Under assumptions (2.5) and (2.6), for all $0 \leq i \leq p$,*

$$\sup_{1 \leq t \leq n/2} \frac{1}{t^{i-\Delta+1/2}} \left| \sum_{j=1}^{\lfloor t \rfloor} j^i \varepsilon_j - \sigma \int_0^t x^i \, dW_{1,n}(x) \right| = \mathcal{O}_P(1), \qquad (4.1)$$

$$\sup_{n/2 \leq t \leq n-1} \frac{1}{(n-t)^{i-\Delta+1/2}} \left| \sum_{j=\lfloor t \rfloor}^n j^i \varepsilon_j - \sigma \int_t^n x^i \, dW_{2,n}(x) \right| = \mathcal{O}_P(1) \qquad (4.2)$$

*as $n \to \infty$. Also, for each $n$,*

$$\left\{ \int_0^s x^i \, dW_{1,n}(x) : s \geq 0 \right\} \stackrel{\mathcal{D}}{=} \left\{ W\left( \frac{s^{2i+1}}{2i+1} \right) : s \geq 0 \right\}, \qquad (4.3)$$

$$\left\{ \int_s^n x^i \, dW_{2,n}(x) : 0 \leq s \leq n-1 \right\} \stackrel{\mathcal{D}}{=} \left\{ W\left( \frac{(n-s)^{2i+1}}{2i+1} \right) : 0 \leq s \leq n-1 \right\}. \qquad (4.4)$$

**Proof.** Note that, by the modulus of continuity of a Brownian motion (see Csörgő and Révész (1981), Theorem 1.21), condition (2.5) can be rewritten as

$$\sup_{1 \leq t \leq n/2} \frac{1}{t^{1/2-\Delta}} \left| \sum_{j=1}^{\lfloor t \rfloor} \varepsilon_j - \sigma W_{1,n}(t) \right| = \mathcal{O}_P(1) \qquad (n \to \infty).$$

Applying integration by parts hence yields (4.1). The proof of (4.2) follows in a similar fashion from the approximation in (2.6).

The statements given in (4.3) and (4.4) follow from computing the corresponding covariance functions. This suffices since the integral processes in (4.1) and (4.2) are Gaussian with mean zero. □

Observe that, by definition, $\mathbf{C}_n = \mathbf{C}_k + \tilde{\mathbf{C}}_k$. Therefore,

$$\mathbf{S}_k^T \tilde{\mathbf{C}}_k^{-1} \mathbf{C}_n \tilde{\mathbf{C}}_k^{-1} \mathbf{S}_k = \mathbf{S}_k^T \tilde{\mathbf{C}}_k^{-1} \mathbf{S}_k + \mathbf{S}_k^T \mathbf{C}_k^{-1} \mathbf{S}_k,$$

with the two terms on the right-hand side being the subject of study in the following. Let $\|\cdot\|$ denote the maximum norm of both vectors and matrices. We obtain the following orders of magnitude for the matrices $\mathbf{C}_k$ and the inverse matrices $\tilde{\mathbf{C}}_k^{-1}$.

**Lemma 4.2.** *As $n \to \infty$,*

$$\max_{1 \leq k \leq n} \left\|\frac{1}{k}\mathbf{C}_k\right\| = \mathcal{O}(1) \quad and \quad \max_{1 \leq k \leq n/2} \|\tilde{\mathbf{C}}_k^{-1}\| = \mathcal{O}\left(\frac{1}{n}\right). \tag{4.5}$$

**Proof.** The first statement in (4.5) follows directly from the definition of $\mathbf{C}_k$. For the second statement, elementary approximations of sums with integrals imply

$$\sup_{0 \leq t \leq 1/2} \left\|\frac{1}{n}\tilde{\mathbf{C}}_{\lfloor nt \rfloor} - \bar{\mathbf{C}}_t\right\| = o(1) \qquad (n \to \infty),$$

where, for $t \in [0, 1/2]$, the matrix $\bar{\mathbf{C}}_t = \{\bar{C}_t(i,j) : 0 \leq i, j \leq p\}$ is defined by

$$\bar{C}_t(i,j) = \int_t^1 x^{i+j}\,\mathrm{d}x = \frac{1}{i+j+1}(1 - t^{i+j+1}).$$

By definition, $\bar{\mathbf{C}}_t$ is continuous on $[0, 1/2]$. Observe, moreover, that $\bar{\mathbf{C}}_t$ is the covariance matrix of the random variables $\int_t^1 x^i\,\mathrm{d}W(x)$, $0 \leq i \leq p$, where $\{W(s) : s \geq 0\}$ denotes a standard Brownian motion. Since these variables are linearly independent, $\bar{\mathbf{C}}_t$ is non-singular for any $t \in [0, 1/2]$. Hence, the proof is complete. □

Next, observe that the vector $\mathbf{S}_k$ can be decomposed into a partial sum vector of $\varepsilon_i$'s, each of the latter random variables weighted with the corresponding polynomial regressor $\mathbf{x}_i$ ($1 \leq i \leq k$) and a second term consisting of $\mathbf{C}_k$ and the difference $\hat{\boldsymbol{\beta}}_n - \boldsymbol{\beta}$, that is,

$$\mathbf{S}_k = \mathbf{v}_k - \mathbf{C}_k(\hat{\boldsymbol{\beta}}_n - \boldsymbol{\beta}), \qquad \mathbf{v}_k = \sum_{i=1}^k \mathbf{x}_i \varepsilon_i. \tag{4.6}$$

To investigate the contributing range of indices $k$, let $\alpha, \beta > 0$ and define

$$a(n) = \log^\alpha n \quad \text{and} \quad b(n) = \frac{n}{\log^\beta n} \qquad (n \geq 1).$$

We obtain the following lemma which shows that the term $\mathbf{S}_k^T \tilde{\mathbf{C}}_k^{-1} \mathbf{S}_k$ is asymptotically negligible.

**Lemma 4.3.** *If (2.5) and (2.6) hold, then, as $n \to \infty$,*

$$\max_{1 \leq k \leq n/2} \mathbf{S}_k^T \tilde{\mathbf{C}}_k^{-1} \mathbf{S}_k = \mathcal{O}_P(1), \tag{4.7}$$

$$\max_{a(n) \leq k \leq b(n)} \mathbf{S}_k^T \tilde{\mathbf{C}}_k^{-1} \mathbf{S}_k = \mathcal{O}_P(1) \log^{-\beta} n. \tag{4.8}$$

**Proof.** By Lemma 4.2 and (4.6), as $n \to \infty$,

$$\max_{1 \leq k \leq n/2} \mathbf{S}_k^T \tilde{\mathbf{C}}_k^{-1} \mathbf{S}_k = \mathcal{O}_P\left(\frac{1}{n}\right) \max_{1 \leq k \leq n/2} \|\mathbf{S}_k\|^2$$

$$= \mathcal{O}_P\left(\frac{1}{n}\right) \left( \max_{1 \leq k \leq n/2} \|\mathbf{v}_k\|^2 + \|\hat{\boldsymbol{\beta}}_n - \boldsymbol{\beta}\|^2 \max_{1 \leq k \leq n/2} \|\mathbf{C}_k\|^2 \right).$$

Now, the central limit theorem applied to the sequence $\{\hat{\boldsymbol{\beta}}_n - \boldsymbol{\beta}\}$ implies that

$$\|\sqrt{n}(\hat{\boldsymbol{\beta}}_n - \boldsymbol{\beta})\| = \mathcal{O}_P(1) \qquad (n \to \infty). \tag{4.9}$$

Hence, by Lemma 4.2,

$$\|\hat{\boldsymbol{\beta}}_n - \boldsymbol{\beta}\|^2 \max_{1 \leq k \leq n/2} \|\mathbf{C}_k\|^2 = \mathcal{O}_P(n) \qquad (n \to \infty).$$

On using Lemma 4.1, we get (see the proof of (4.8) below and proceed in the same way)

$$\max_{1 \leq k \leq n/2} \|\mathbf{v}_k\|^2 = \mathcal{O}_P(n) \qquad (n \to \infty),$$

completing the proof of (4.7).

Similarly, by Lemma 4.2,

$$\max_{1 \leq k \leq b(n)} \mathbf{S}_k^T \tilde{\mathbf{C}}_k^{-1} \mathbf{S}_k = \mathcal{O}_P\left(\frac{1}{n}\right) \left( \max_{1 \leq k \leq b(n)} \|\mathbf{v}_k\|^2 + \mathcal{O}_P\left(\frac{1}{n}\right) \max_{1 \leq k \leq b(n)} \|\mathbf{C}_k\|^2 \right)$$

$$= \mathcal{O}_P\left(\frac{1}{n}\right) \max_{1 \leq k \leq b(n)} \|\mathbf{v}_k\|^2 + \mathcal{O}_P(\log^{-2\beta} n).$$

Applying Lemma 4.1 yields that

$$\max_{1 \leq k \leq b(n)} \|\mathbf{v}_k\| = \sigma \max_{0 \leq i \leq p} \max_{1 \leq k \leq b(n)} \frac{1}{n^i} \left| \int_0^k x^i \, dW_{1,n}(x) \right|$$

$$+ \mathcal{O}_P(1) \max_{0 \leq i \leq p} \max_{1 \leq k \leq b(n)} n^{-i} (b(n))^{1/2 + i - \Delta},$$

where, inserting the definition of $b(n)$, the second term can be estimated by

$$\max_{0 \leq i \leq p} \frac{1}{n^i} \left( \frac{n}{\log^\beta n} \right)^{1/2 + i - \Delta} \leq \left( \frac{n}{\log^\beta n} \right)^{1/2 - \Delta}.$$

Since, for all $0 \leq i \leq p$,

$$\frac{1}{(b(n))^{i+1/2}} \max_{0 \leq t \leq b(n)} \left| \int_0^t x^i \, dW(x) \right| \overset{\mathcal{D}}{=} \frac{1}{\sqrt{2i+1}} \sup_{0 \leq t \leq 1} \left| \int_0^t x^i \, dW(x) \right|,$$

we get that

$$\max_{0\leq i\leq p}\max_{0\leq t\leq b(n)}\frac{1}{n^i}\left|\int_0^t x^i\,dW_{1,n}(x)\right|$$
$$=\mathcal{O}_P(1)\max_{0\leq i\leq p}\max_{0\leq t\leq b(n)}\frac{1}{n^i}\left(\frac{n}{\log^\beta n}\right)^{i+1/2}=\mathcal{O}_P(1)\left(\frac{n}{\log^\beta n}\right)^{1/2},$$

completing the proof of (4.8). □

It remains to examine the term $\mathbf{S}_k^T \mathbf{C}_k^{-1} \mathbf{S}_k$. Applying the expression for $\mathbf{S}_k$ obtained in (4.6) leads to the refined decomposition

$$\mathbf{S}_k^T \mathbf{C}_k^{-1} \mathbf{S}_k = (\mathbf{v}_k^T - (\hat{\boldsymbol{\beta}}_n - \boldsymbol{\beta})^T \mathbf{C}_k)\mathbf{C}_k^{-1}(\mathbf{v}_k - \mathbf{C}_k(\hat{\boldsymbol{\beta}}_n - \boldsymbol{\beta}))$$
$$= \mathbf{v}_k^T \mathbf{C}_k^{-1}\mathbf{v}_k - 2(\hat{\boldsymbol{\beta}}_n - \boldsymbol{\beta})^T \mathbf{v}_k + (\hat{\boldsymbol{\beta}}_n - \boldsymbol{\beta})^T \mathbf{C}_k (\hat{\boldsymbol{\beta}}_n - \boldsymbol{\beta}).$$

Orders of magnitude of quantities appearing on the right-hand side of the latter equation array are provided in the following lemma. We restrict the discussion to the range $k\leq n/2$, the other half of the observations, for which $k>n/2$, will be dealt with after Lemma 4.12, repeating the arguments developed now.

**Lemma 4.4.** *If (2.5) and (2.6) hold, then, as $n\to\infty$,*

$$\max_{1\leq k\leq n/2}(\hat{\boldsymbol{\beta}}_n - \boldsymbol{\beta})^T \mathbf{C}_k (\hat{\boldsymbol{\beta}}_n - \boldsymbol{\beta}) = \mathcal{O}_P(1), \tag{4.10}$$

$$\max_{1\leq k\leq b(n)}(\hat{\boldsymbol{\beta}}_n - \boldsymbol{\beta})^T \mathbf{C}_k (\hat{\boldsymbol{\beta}}_n - \boldsymbol{\beta}) = \mathcal{O}_P(1)\log^{-\beta} n, \tag{4.11}$$

$$\max_{1\leq k\leq n/2}(\hat{\boldsymbol{\beta}}_n - \boldsymbol{\beta})^T \mathbf{v}_k = \mathcal{O}_P(1), \tag{4.12}$$

$$\max_{1\leq k\leq b(n)}(\hat{\boldsymbol{\beta}}_n - \boldsymbol{\beta})^T \mathbf{v}_k = \mathcal{O}_P(1)\log^{-\beta/2} n. \tag{4.13}$$

**Proof.** Using Lemma 4.2 and (4.9), we conclude that statement (4.10) holds. For (4.11), it can be similarly obtained that

$$\max_{1\leq k\leq b(n)}(\hat{\boldsymbol{\beta}}_n - \boldsymbol{\beta})^T \mathbf{C}_k(\hat{\boldsymbol{\beta}}_n - \boldsymbol{\beta}) = \mathcal{O}_P\left(\frac{b(n)}{n}\right)\qquad (n\to\infty).$$

Lemma 4.2 implies that, as $n\to\infty$,

$$\frac{1}{\sqrt{n}}\max_{1\leq k\leq n/2}\|\mathbf{v}_k\| = \mathcal{O}_P(1) \quad\text{and}\quad \frac{1}{\sqrt{b(n)}}\max_{1\leq k\leq b(n)}\|\mathbf{v}_k\| = \mathcal{O}_P(1).$$

Therefore, statements (4.12) and (4.13) follow from Lemma 4.1 and (4.9). □

Next, we consider the maximum of $\mathbf{S}_k^T \mathbf{C}_k^{-1} \mathbf{S}_k$ on three different ranges determined by $a(n)$ and $b(n)$. To this end, we write

$$\max_{1 \leq k \leq n/2} \mathbf{S}_k^T \mathbf{C}_k^{-1} \mathbf{S}_k = \max\{z_{n,1}, z_{n,2}, z_{n,3}\}$$

with

$$z_{n,1} = \max_{1 \leq k \leq a(n)} \mathbf{S}_k^T \mathbf{C}_k^{-1} \mathbf{S}_k, \qquad z_{n,2} = \max_{a(n) \leq k \leq b(n)} \mathbf{S}_k^T \mathbf{C}_k^{-1} \mathbf{S}_k$$

and

$$z_{n,3} = \max_{b(n) \leq k \leq n/2} \mathbf{S}_k^T \mathbf{C}_k^{-1} \mathbf{S}_k.$$

Lemma 4.5 provides approximations of $z_{n,1}$, $z_{n,2}$ and $z_{n,3}$ using the partial sum vectors $\mathbf{v}_k$ introduced in (4.6). It turns out that the maximum on all three ranges is completely determined by the maximum of the quadratic form $\mathbf{v}_k^T \mathbf{C}_k^{-1} \mathbf{v}_k$.

**Lemma 4.5.** *If (2.5) and (2.6) hold, then, as $n \to \infty$,*

$$z_{n,1} = \max_{1 \leq k \leq a(n)} \mathbf{v}_k^T \mathbf{C}_k^{-1} \mathbf{v}_k + \mathcal{O}_P(1),$$

$$z_{n,2} = \max_{a(n) \leq k \leq b(n)} \mathbf{v}_k^T \mathbf{C}_k^{-1} \mathbf{v}_k + \mathcal{O}_P(\log^{-\beta/2} n),$$

$$z_{n,3} = \max_{b(n) \leq k \leq n/2} \mathbf{v}_k^T \mathbf{C}_k^{-1} \mathbf{v}_k + \mathcal{O}_P(1).$$

**Proof.** All three approximations follow immediately from the foregoing Lemma 4.4. □

Introducing the notation

$$\begin{aligned}
\mathbf{s}_k^T &= (s_{k,0}, \ldots, s_{k,p}), \qquad s_{k,i} = \sum_{j=1}^k j^i \varepsilon_j, \\
\mathbf{B}_n &= \{B_n(i,j) : 0 \leq i, j \leq p\}, \\
B_n(i,i) &= n^{-i}, \qquad B_n(i,j) = 0 \qquad (i \neq j), \\
\mathbf{D}_k^{-1} &= \mathbf{B}_n \mathbf{C}_k^{-1} \mathbf{B}_n,
\end{aligned} \qquad (4.14)$$

we obtain, for all $1 \leq k < n$,

$$\mathbf{v}_k^T \mathbf{C}_k^{-1} \mathbf{v}_k = \mathbf{s}_k^T \mathbf{B}_n \mathbf{C}_k^{-1} \mathbf{B}_n \mathbf{s}_k = \mathbf{s}_k^T \mathbf{D}_k^{-1} \mathbf{s}_k.$$

The diagonal matrix $\mathbf{B}_n$ simplifies the partial sums $\mathbf{v}_k$ by taking out the norming factor $n^{-i}$, returning the partial sums $\mathbf{s}_k$ which are used to rewrite the quadratic form, as done in the last display above. Note that the vectors $\mathbf{s}_k$ possess covariance matrix $\sigma^2 \mathbf{D}_k$.

We now transform $\mathbf{s}_k$ into a new vector $\tilde{\mathbf{s}}_k$ which has uncorrelated components via an application of the Gram–Schmidt orthonormalization. To define $\tilde{\mathbf{s}}_k = (\tilde{s}_{k,0}, \ldots, \tilde{s}_{k,p})^T$, let $\tilde{s}_{k,0} = s_{k,0}$ and set, recursively, for $1 \leq i \leq p$,

$$\tilde{s}_{k,i} = s_{k,i} - \alpha_{k,i,0}\tilde{s}_{k,0} - \cdots - \alpha_{k,i,i-1}\tilde{s}_{k,i-1}. \tag{4.15}$$

The coefficients $\alpha_{k,i,j}$ are given as projections of $s_{k,i}$ onto the new variables $\tilde{s}_{k,j}$ ($j < i$), that is,

$$\alpha_{k,i,j} = \frac{\mathrm{E}s_{k,i}\tilde{s}_{k,j}}{\sqrt{\mathrm{E}\tilde{s}_{k,j}^2}}, \qquad 0 \leq j < i \leq p.$$

Consequently, $\mathrm{E}\tilde{s}_{k,1}\tilde{s}_{k,j} = 0$ for $0 \leq i \neq j \leq p$. According to Rao (1973), the matrices $\sigma^{-1}\mathbf{D}_k^{-1/2}$ are related to the Gram–Schmidt orthonormalization in (4.15) via the equality

$$\mathbf{s}_k^T \sigma^{-1}\mathbf{D}_k^{-1/2} = \left(\frac{\tilde{s}_{k,0}}{\sqrt{\mathrm{Var}\,\tilde{s}_{k,0}}}, \ldots, \frac{\tilde{s}_{k,p}}{\sqrt{\mathrm{Var}\,\tilde{s}_{k,p}}}\right). \tag{4.16}$$

Next, we study moments of $s_{k,i}$ and $\tilde{s}_{k,i}$. Note that the precise limits given in the following lemma will be specified in Lemma 4.10 below.

**Lemma 4.6.** *If assumption (2.4) is satisfied, then, for all $0 \leq i \leq p$,*

$$\lim_{k\to\infty} \frac{\mathrm{E}s_{k,i}^2}{k^{2i+1}} = \frac{\sigma^2}{2i+1}, \tag{4.17}$$

$$\lim_{k\to\infty} \frac{\mathrm{E}\tilde{s}_{k,i}^2}{k^{i+1}} = \beta_i > 0. \tag{4.18}$$

*Furthermore, for all $0 \leq j < i \leq p$,*

$$\lim_{k\to\infty} \frac{\mathrm{E}s_{k,i}\tilde{s}_{k,j}}{k^{i+j+1}} = \beta_{i,j}, \tag{4.19}$$

$$\lim_{k\to\infty} \frac{\alpha_{k,i,j}}{k^{i-j}} = \alpha_{i,j}. \tag{4.20}$$

**Proof.** Since, by assumption (2.4) on the errors $\{\varepsilon_i\}$,

$$\mathrm{E}s_{k,i}^2 = \sigma^2 \sum_{\ell=1}^{k} \ell^{2i},$$

(4.17) follows after elementary calculations. Relations (4.18)–(4.20) can be established using mathematical induction and the definition of the coefficients $\alpha_{k,i,j}$. Since the precise limits are obtained below, we only show here that, in (4.18), the limit $\beta_i > 0$. To this

end, writing

$$\tilde{s}_{k,i} = \sum_{\ell=1}^{k} (\ell^i - \alpha_{k,i,0} - \alpha_{k,i,1}\ell - \cdots - \alpha_{k,i,i-1}\ell^{i-1})\varepsilon_\ell,$$

assumption (2.4) yields that

$$\mathrm{E}\tilde{s}_{k,i}^2 = \sigma^2 \sum_{\ell=1}^{k} (\ell^i - \alpha_{k,i,0} - \alpha_{k,i,1}\ell - \cdots - \alpha_{k,i,i-1}\ell^{i-1})^2.$$

Therefore,

$$\lim_{k\to\infty} \frac{\mathrm{E}\tilde{s}_{k,i}^2}{k^{2i+1}} = \int_0^1 (t^i - \alpha_{i,0} - \alpha_{i,1}t - \cdots - \alpha_{i,i-1}t^{i-1})^2 \, \mathrm{d}t$$

so that $\beta_i > 0$, since the integrand is a non-negative, non-constant polynomial. $\square$

**Lemma 4.7.** *If (2.4)–(2.6) hold, then, as $n \to \infty$,*

$$\max_{1 \le k \le a(n)} \mathbf{s}_k^T \mathbf{D}_k^{-1} \mathbf{s}_k = \mathcal{O}_P(\log\log\log n), \tag{4.21}$$

$$\max_{b(n) \le k \le n/2} \mathbf{s}_k^T \mathbf{D}_k^{-1} \mathbf{s}_k = \mathcal{O}_P(\log\log\log n) \tag{4.22}$$

*and*

$$z_{n,1} = \mathcal{O}_P(\log\log\log n), \tag{4.23}$$

$$z_{n,3} = \mathcal{O}_P(\log\log\log n). \tag{4.24}$$

**Proof.** According to Lemma 4.5, it is enough to prove (4.21) and (4.22). Using Lemma 4.1 and the Darling–Erdős law for Brownian motions (cf. Csörgő and Révész (1981)), we get, for all $0 \le i \le p$,

$$\max_{1 \le k \le a(n)} \frac{1}{k^{i+1/2}} \left| \sum_{\ell=1}^{k} \ell^i \varepsilon_\ell \right| = \mathcal{O}_P(\sqrt{\log\log a(n)}) \qquad (n \to \infty).$$

Hence, (4.22) follows from (4.15), (4.16) and Lemma 4.6. $\square$

(4.23) and (4.24) imply that, for any $c > 0$,

$$\max\{z_{n,1}, z_{n,3}\} - c\log\log n \xrightarrow{P} -\infty \qquad (n \to \infty)$$

so that it suffices to consider the asymptotics for the remaining term $z_{n,2}$. In what follows, we are going to approximate $\mathbf{s}_k^T \mathbf{D}_k^{-1} \mathbf{s}_k$ with a quadratic form consisting of normal random vectors. To this end, set $N_i = N_{i,n} = W_{1,n}(i) - W_{1,n}(i-1)$, where $\{W_{n,1}(s) : s \geq 0\}$ is the Brownian motion defined in assumption (2.5). Let

$$\mathbf{n}_k^T = (n_{k,0}, \ldots, n_{k,p}), \qquad n_{k,i} = \sum_{\ell=1}^{k} \ell^i N_i.$$

**Lemma 4.8.** *If (2.4)–(2.6) hold, then, as $n \to \infty$,*

$$\frac{1}{\sigma^2} \max_{a(n) \leq k \leq b(n)} \mathbf{s}_k^T \mathbf{D}_k^{-1} \mathbf{s}_k = \max_{a(n) \leq k \leq b(n)} \mathbf{n}_k^T \mathbf{D}_k^{-1} \mathbf{n}_k + \mathcal{O}_P(\log^{-\alpha\Delta} n \sqrt{\log \log n}).$$

**Proof.** For $k \geq 1$, set $\tilde{n}_{k,0} = n_{k,0}$ and define, recursively, for $1 \leq i \leq p$,

$$\tilde{n}_{k,i} = n_{k,i} - \alpha_{k,i,0} \tilde{n}_{k,0} - \cdots - \alpha_{k,i,i-1} \tilde{n}_{k,i-1}. \tag{4.25}$$

Following the arguments leading to (4.15) and (4.16), we get that

$$\mathbf{n}_k^T \mathbf{D}_k^{-1/2} = \left( \frac{\tilde{n}_{k,0}}{\sqrt{\mathrm{E} \tilde{n}_{k,0}^2}}, \ldots, \frac{\tilde{n}_{k,p}}{\sqrt{\mathrm{E} \tilde{n}_{k,p}^2}} \right). \tag{4.26}$$

From the proof of Lemma 4.1, we obtain, as $n \to \infty$,

$$\max_{1 \leq k \leq n/2} \frac{1}{k^{i+1/2-\Delta}} \left| \frac{1}{\sigma} \sum_{\ell=1}^{k} \ell^i \varepsilon_\ell - \sum_{\ell=1}^{k} \ell^i N_i \right| = \mathcal{O}_P(1). \tag{4.27}$$

Combining (4.15), (4.16) and (4.26) with Lemma 4.6 leads to

$$\max_{a(n) \leq k \leq n/2} \|(\sigma^{-1} \mathbf{s}_k^T - \mathbf{n}_k^T) \mathbf{D}_k^{-1/2}\| = \mathcal{O}_P(\log^{-\alpha\Delta} n).$$

By the Darling–Erdős law for the Wiener process and Lemma 4.1, we get

$$\max_{1 \leq k \leq n/2} \|\sigma^{-1} \mathbf{s}_k^T \mathbf{D}_k^{-1/2}\| = \mathcal{O}_P(\sqrt{\log \log n}) = \max_{1 \leq k \leq n/2} \|\sigma^{-1} \mathbf{n}_k^T \mathbf{D}_k^{-1/2}\|.$$

Therefore, the lemma is proved. □

Next, we replace the discrete normal partial sums $n_{k,i}$ with continuous integrals. Define

$$\mathbf{q}_n^T(t) = (q_{n,0}(t), \ldots, q_{n,p}(t)), \qquad q_{n,i}(t) = \int_0^t x^i \, \mathrm{d}W_{1,n}(x),$$

with $\{W_{1,n}(s) : s \geq 0\}$ being the Brownian motion defined in (2.5).

**Lemma 4.9.** *If (2.5) and (2.6) hold, then, as $k \to \infty$,*

$$\max_{k \leq t < k+1} |q_{n,i}(t) - n_{k,i}| = \mathcal{O}_P(k^i \sqrt{\log k}).$$

**Proof.** Using the modulus of continuity of a Brownian motion (cf. Csörgő and Révész (1981), Theorem 1.2.1) and integration by parts, we obtain that

$$\max_{k \leq t < k+1} \left| \int_0^t x^i \, \mathrm{d}(W_{1,n}(x) - W_{1,n}(\lfloor x \rfloor)) \right|$$

$$\leq \sup_{k \leq t < k+1} t^i |W_{1,n}(t) - W_{1,n}(\lfloor t \rfloor)| + \left| \int_0^1 x^i \, \mathrm{d}W_{1,n}(x) \right|$$

$$+ \sup_{k \leq t < k+1} \left| \int_0^t i x^{i-1} (W_{1,n}(x) - W_{1,n}(\lfloor x \rfloor)) \, \mathrm{d}x \right|$$

$$= \mathcal{O}_P(k^i \sqrt{\log k}) + \mathcal{O}_P(1) \int_0^k x^{i-1} \sqrt{\log(|x|+2)} \, \mathrm{d}x$$

$$= \mathcal{O}_P(k^i \sqrt{\log k}),$$

proving the result. $\square$

Copying the arguments leading to (4.15), we define $\tilde{q}_{n,0}(t) = q_{n,0}(t)$ and, recursively, for $1 \leq i \leq p$,

$$\tilde{q}_{n,i}(t) = q_{n,i} - \alpha_{i,0}(t)\tilde{q}_{n,0} - \cdots - \alpha_{i,i-1}(t)\tilde{q}_{i,i-1}(t).$$

Observe that since the distribution of $\{W_{1,n}(s) : s \geq 0\}$ is independent of $n$, so are the coeffcients $\alpha_{i,j}(t)$. By construction, $\mathrm{E}\tilde{q}_{n,i}(t)\tilde{q}_{n,j}(t) = 0$ for $0 \leq i \neq j \leq p$ because, clearly,

$$\alpha_{i,j}(t) = \frac{\mathrm{E} q_{n,i}(t)\tilde{q}_{n,j}(t)}{\mathrm{E}\tilde{q}_{n,i}^2(t)}, \qquad 0 \leq j < i \leq p.$$

Denote by

$$\bar{\mathbf{D}}_t = \{\bar{D}_t(i,j) : 0 \leq i, j \leq p\}, \qquad \bar{D}_t(i,j) = \int_0^t x^{i+j} \, \mathrm{d}x,$$

the covariance matrix of the random vector $\mathbf{q}_n(t)$.

**Lemma 4.10.** *If (2.4) holds, then*

$$\alpha_{i,j}(t) = t^{i-j}\alpha_{i,j}, \qquad 0 \leq j < i \leq p,$$

*where $\alpha_{i,j}$ are constants.*

**Proof.** This is easily verified using mathematical induction or the scale transformation of Brownian motions. □

**Lemma 4.11.** *If (2.4)–(2.6) hold, then, as $n \to \infty$,*

$$\sup_{a(n) \leq t \leq n/2} |\mathbf{n}_{\lfloor t \rfloor}^T \mathbf{D}_{\lfloor t \rfloor}^{-1/2} - \mathbf{q}_n^T(t) \mathbf{D}_t^{-1/2}| = \mathcal{O}_P(\log^{-\alpha/2} \sqrt{\log \log n}).$$

**Proof.** Recall (4.26) and note that the difference between the coefficients $\alpha_{k,i,j}$ and $\alpha_{i,j}(t)$ is that the first expressions are defined as Riemann sums which approximate the integrals defining the latter. It can, hence, be easily verified that

$$\sup_{k \leq t < k+1} \frac{|\alpha_{i,j}(t) - \alpha_{k,i,j}|}{\alpha_{k,i,j}} = \mathcal{O}\left(\frac{1}{k}\right) \qquad (k \to \infty).$$

The result now follows from Lemmas 4.8–4.10. □

Define

$$\hat{\mathbf{q}}_n(t) = (\hat{q}_0(t), \ldots, \hat{q}_p(t)), \qquad \hat{q}_i(t) = \frac{\int_0^t g_{i,t}(x) \, dW_{1,n}(x)}{(\mathrm{E}[\int_0^t g_{i,t}(x) \, dW_{1,n}(x)]^2)^{1/2}},$$

where $\{W_{1,n}(s) : s \geq 0\}$ is defined in (2.5) and

$$g_{i,t}(x) = x^i - \alpha_{i,0} t^i - \alpha_{i,1} t^{i-1} x - \cdots - \alpha_{i,i-1} t x^{i-1}.$$

Note that the distribution of $\hat{\mathbf{q}}_n(t)$ does not depend on $n$. Thus, we arrive at the following result.

**Lemma 4.12.** *If (2.4)–(2.6) hold, then, as $n \to \infty$,*

$$\sigma^{-2} z_{n,2} = \sup_{a(n) \leq t \leq b(n)} \hat{\mathbf{q}}_n^T(t) \hat{\mathbf{q}}_n(t) + \mathcal{O}_P(\sqrt{\log \log n}[\log^{-\alpha \Delta} n + \log n^{-\alpha/2} n]).$$

**Proof.** This follows on combining the results of the foregoing lemmas. □

On summarizing what has been proven thus we conclude that, restricted to $1 \leq k \leq n/2$, the test statistic $T_n$ introduced in (2.2) can be asymptotically characterized by the supremum (over a restricted range determined by $a(n)$ and $b(n)$) of the inner product $\hat{\mathbf{q}}_n^T(t) \hat{\mathbf{q}}_n(t)$, where the vectors $\hat{\mathbf{q}}_n(t)$ have components given by the (normalized) integrals $\int_0^t g_{i,t}(x) \, dW_{1,n}(x)$. Clearly, the functions $g_{i,t}(x)$ are polynomials of order $i$, whose leading coefficient (i.e., the coefficient assigned to the monomial $x^i$) is 1 and, by construction, they possess the property $\int_0^t g_{i,t}(x) g_{j,t}(x) \, dx = 0$ if $i \neq j$. So, the functions $g_{i,t}(x)$ are exactly the Legendre polynomials on the interval $[0, t]$ with normalized coefficient for the monomial $x^i$. This completes the proof for the range $1 \leq k \leq n/2$.

In the remainder of the subsection, we consider $\max_{n/2<k<n} \mathbf{S}_k^T \mathbf{C}_k^{-1} \mathbf{S}_k$, adapting the arguments developed before. First, observe that

$$\mathbf{S}_k = \tilde{\mathbf{C}}_k \mathbf{C}_n^{-1} \sum_{i=1}^{n} \mathbf{x}_i y_i - \sum_{i=k+1}^{n} \mathbf{x}_i y_i.$$

Let

$$c(n) = n - \frac{n}{\log^\beta n} \quad \text{and} \quad d(n) = n - \log^\alpha n.$$

Along the lines of the proofs of Lemmas 4.1 and 4.2–4.7, one can establish the following orders of magnitude:

$$\max_{n/2 \le k \le c(n)} \mathbf{S}_k^T \mathbf{C}_k^{-1} \mathbf{S}_k = \mathcal{O}_P(\log \log \log n);$$

$$\max_{d(n) \le k < n-p} \mathbf{S}_k^T \mathbf{C}_k^{-1} \mathbf{S}_k = \mathcal{O}_P(\log \log \log n).$$

Let

$$\mathbf{r}_k^T = (r_{k,0}, \ldots, r_{k,p}), \qquad r_{k,i} = \sum_{j=k+1}^{n-1} j^i \varepsilon_j$$

and $\mathbf{R}_k^{-1} = \mathbf{B}_n \tilde{\mathbf{C}}_k^{-1} \mathbf{B}_n$, with $\mathbf{B}_n$ defined in (4.14). As in Lemma 4.5, we obtain that

$$\max_{c(n) \le k \le d(n)} \mathbf{S}_k^T \mathbf{C}_k^{-1} \mathbf{S}_k = \max_{c(n) \le k \le d(n)} \mathbf{r}_k^T \mathbf{R}_k^{-1} \mathbf{r}_k + \mathcal{O}_P(\log^{-\beta/2} n).$$

Repeating the arguments that lead to the statements contained in Lemmas 4.6–4.12 here yields

$$\frac{1}{\sigma^2} \max_{c(n) \le k \le d(n)} \mathbf{r}_k^T \mathbf{R}_k^{-1} \mathbf{r}_k = \sup_{c(n) \le t \le d(n)} \mathbf{e}^T(t) \hat{\mathbf{R}}_t^{-1} \mathbf{e}(t),$$

where, using the Brownian motion $\{W_{2,n}(s) : s \ge 0\}$ specified in assumption (2.6),

$$\mathbf{e}^T(t) = (e_0(t), \ldots, e_p(t)), \qquad e_i(t) = \int_t^n x^i \, dW_{2,n}(x), \qquad (4.28)$$

$$\hat{\mathbf{R}}_t = \{\hat{R}_t(i,j) : 0 \le i, j \le p\}, \qquad \hat{R}_t(i,j) = \int_t^n x^{i+j} \, dx.$$

If we finally define

$$\tilde{\mathbf{e}}^T(t) = (\tilde{e}_0(t), \ldots, \tilde{e}_p(t)), \qquad \tilde{e}_i(t) = \int_0^t (n-x)^i \, dW(x),$$

$$\tilde{\mathbf{R}}_t = \{\hat{R}_t(i,j) : 0 \le i, j \le p\}, \qquad \tilde{R}_t(i,j) = \int_0^t (n-x)^{i+j} \, dx,$$

where $\{W(s): s \geq 0\}$ denotes a standard Brownian motion, then

$$\sup_{c(n) \leq t \leq d(n)} \mathbf{e}^T(t)\hat{\mathbf{R}}_t^{-1}\mathbf{e}(t) \stackrel{\mathcal{D}}{=} \sup_{a(n) \leq t \leq b(n)} \tilde{\mathbf{e}}^T(t)\tilde{\mathbf{R}}_t^{-1}\tilde{\mathbf{e}}(t).$$

The orthonormalization in $\tilde{\mathbf{e}}^T(t)\tilde{\mathbf{R}}_t^{-1/2}$ starts with $W(t)$ and the coordinates are uncorrelated, the $i$th coordinate being the integral of an $i$th order polynomial with respect to $\{W(s): s \geq 0\}$. Hence, we are again on the restricted range determined by $c(n)$ and $d(n)$, integrating the Legendre polynomials on $[0, t]$ with respect to the Brownian motion $\{W(s): s \geq 0\}$.

The final lemma establishes that the suprema taken over the first part, $a(n) < k < b(n)$, and second part, $c(n) < k < d(n)$, have the same distribution, but are independent.

**Lemma 4.13.** *For each $n$, there are $\{\tilde{\mathbf{e}}_n(t), 0 \leq t \leq n\} = \{\tilde{\mathbf{e}}(t), 0 \leq t \leq n\}$ such that*

$$\sup_{a(n) \leq t \leq b(n)} \tilde{\mathbf{e}}_n^T(t)\tilde{\mathbf{R}}_t^{-1}\tilde{\mathbf{e}}_n(t) \stackrel{\mathcal{D}}{=} \sup_{a(n) \leq t \leq b(n)} \hat{\mathbf{q}}_n^T(t)\hat{\mathbf{q}}_n(t),$$

*where the first and second suprema are independent. Moreover, as $n \to \infty$,*

$$\max_{c(n) \leq k \leq d(n)} \mathbf{S}_k^T \mathbf{C}_k^{-1} \mathbf{S}_k = \sup_{a(n) \leq t \leq b(n)} \tilde{\mathbf{e}}_n^T(t)\tilde{\mathbf{R}}_t^{-1}\tilde{\mathbf{e}}_n(t) + \mathcal{O}_P((\log n)^{-\beta/2}) + \mathcal{O}_P(\sqrt{\log \log n}[\log^{-\alpha\Delta} n + \log^{-\alpha} n]).$$

**Proof.** This is immediate. □

## 5. Proofs of Theorem 2.2 and Corollary 2.1

The final section completes the proof of Theorem 2.2 and its corollary by combining the results obtained in the previous section with the extreme value theory for stochastic integrals of Legendre polynomials with respect to Brownian motions as developed in Aue et al. (2007).

**Proof of Theorem 2.2.** Let

$$c(n, x) = x + 2 \log \log n + (p+1) \log \log \log n - 2 \log\left(\frac{2^{(p+1)/2}\Gamma((p+1)/2)}{p+1}\right),$$

where $\Gamma(t)$ denotes the Gamma function. If follows from Lemmas 4.3–4.13 that, as $n \to \infty$,

$$|P\{T_n \leq c(n, x)\} - P\{\max\{\xi_{n,1}, \xi_{n,2}\} \leq c(n, x)\}| \to 0$$

for all $x$, where $\xi_{n,1}$ and $\xi_{n,2}$ are independent, identically distributed random variables satisfying

$$\xi_{n,1} \stackrel{\mathcal{D}}{=} \xi_{n,2} \stackrel{\mathcal{D}}{=} \sup_{a(n) \leq t \leq b(n)} \hat{q}_n^T(t)\hat{q}_n(t).$$

Therein, the quantities $\hat{q}_n(t)$ are the Gaussian processes introduced in Lemma 4.12. An application of the main result in Aue *et al.* (2007) implies that

$$\lim_{n\to\infty} P\{\xi_{n,i} \le c(n,x)\} = \exp(-e^{-x/2}), \qquad i=1,2,$$

and, consequently, the assertion of Theorem 2.2 follows readily. □

**Proof of Corollary 2.1.** We need only observe that

$$\log\log n - \log\log h(n) = \log\left(\left[1 + \gamma\frac{\log\log n}{\log n}\right]^{-1}\right) \to 0 \qquad (n\to\infty)$$

and

$$\log\log\log n - \log\log\log h(n) = \log\left(1 + \frac{\log\log n - \log\log h(n)}{\log\log h(n)}\right) \to 0$$

since $\log\log n - \log\log h(n) \to 0$ and $\log\log h(n) \to \infty$ as $n\to\infty$. □

# Acknowledgements

This research was supported in part by NSF grants DMS 0413653, DMS 0604670 and DMS 0652420, and grants RGC–HKUST6428/06H, MSM 0021620839, GACR 201/06/0186.